\documentclass{amsart}
\usepackage{amsthm,amssymb,amsmath}
\usepackage{mathtools}
\usepackage[margin=1.2in]{geometry}
\usepackage{xcolor}
\definecolor{darkblue}{rgb}{0,0,0.6}
\usepackage[breaklinks, pdftex, ocgcolorlinks, citecolor=darkblue, filecolor=darkblue, linkcolor=darkblue, urlcolor=black]{hyperref}
\usepackage[capitalize]{cleveref}
\usepackage{enumerate}
\usepackage{color}
\usepackage[matrix,arrow,curve]{xy}
\usepackage{tikz-cd}

\newcommand{\Z}{\mathbb{Z}}

\numberwithin{equation}{section}
\newcounter{commentcounter}

\newtheorem{thm}[equation]{Theorem}
\newtheorem{lem}[equation]{Lemma}
\newtheorem{prop}[equation]{Proposition}
\newtheorem*{claim}{Claim}

\newtheorem{cor}[equation]{Corollary}
\newtheorem{question}[equation]{Question}
\theoremstyle{definition}
\newtheorem{defi}[equation]{Definition}
\theoremstyle{remark}
\newtheorem{rem}[equation]{Remark}
\newtheorem{example}[equation]{Example}

\newcommand{\wt}{\widetilde}
\newcommand{\wh}{\widehat}
\DeclareMathOperator{\coker}{coker}

\DeclareMathOperator{\Hom}{Hom}
\DeclareMathOperator{\id}{Id}
\DeclareMathOperator{\im}{im}

\usepackage{todonotes}

\title[4-manifolds with finite abelian fundamental group]{Homotopy classification of 4-manifolds with finite abelian 2-generator fundamental groups}

\author{Daniel Kasprowski}
\address{School of Mathematical Sciences, University of Southampton, Southampton SO17 1BJ, United Kingdom}
\email{\href{mailto:d.kasprowski@soton.ac.uk}{d.kasprowski@soton.ac.uk}}

\author{Mark Powell}
\address{School of Mathematics and Statistics, University of Glasgow, United Kingdom}
\email{\href{mailto:mark.powell@glasgow.ac.uk}{mark.powell@glasgow.ac.uk}}

\author{Benjamin Ruppik}
\address{Faculty of Mathematics and Natural Sciences, Heinrich Heine University D\"{u}sseldorf, Germany}
\email{\href{mailto:benjamin.ruppik@hhu.de}{benjamin.ruppik@hhu.de}}

\def\subjclassname{\textup{2020} Mathematics Subject Classification}
\expandafter\let\csname subjclassname@1991\endcsname=\subjclassname
\expandafter\let\csname subjclassname@2000\endcsname=\subjclassname

\begin{document}

\begin{abstract}
	We show that for an oriented $4$-dimensional Poincar\'e complex $X$ with finite fundamental group, whose $2$-Sylow subgroup is abelian with at most $2$ generators, the homotopy type of $X$ is determined by its quadratic $2$-type.
\end{abstract}

\maketitle

\section{Introduction}

An oriented 4-dimensional  Poincar\'{e} complex is a finite CW complex with a fundamental class $[X] \in H_4(X;\Z)$ such that
\[- \cap [X] \colon C^{4-*}(X;\Z[\pi_1(X)]) \to C_*(X;\Z[\pi_1(X)])\]
is a chain equivalence.
For many fundamental groups, Hambleton and Kreck \cite[Theorem~1.1]{HK} classified oriented 4-dimensional Poincar\'{e} complexes up to orientation-preserving (o.p.) homotopy equivalence, in terms of isomorphism classes of pairs consisting of the Postnikov 2-type $T_X$ and a homology class  $t \in H_4(T_X)$, the image of $[X]$ under some 3-connected map $X \to T_X$. Baues and Bleile~\cite[Corollary~3.2]{Baues-Bleile} showed that this result indeed holds for all fundamental groups.

One could be tempted to view this as a complete classification. However, it is often desirable in applications to have a classification in terms of more readily computable invariants.  The standard suite of invariants of an oriented $4$-dimensional Poincar\'e complex~$X$ are those which comprise the \emph{quadratic $2$-type}:
\[[\pi_1(X),\pi_2(X),k_X,\lambda_X],\]
where $\pi_2(X)$ is considered as a $\Z[\pi_1(X)]$-module, $k_X\in H^3(\pi_1(X);\pi_2(X))$ is the $k$-invariant 
of $X$ and $\lambda_X$ is the equivariant intersection form on $\pi_2(X)$. Since recording the first three is equivalent to knowing the Postnikov 2-type, investigations into the homotopy classification reduce to asking whether $t \in H_4(T_X)$ is determined by $\lambda_X$.

We recall the history of the homotopy classification problem for oriented simply-connected $4$-dimensional Poincar\'e complexes.
Whitehead and Milnor \cite{whitehead-simply-connected-4mfds,milnor-simply-connected-4mfds} classified them in the simply-connected case, and later Wall~\cite[Theorem~5.4]{Wall-Poincare-cxs} did the same for those with cyclic fundamental groups of prime order.
Hambleton and Kreck \cite{HK} showed that oriented $4$-dimensional Poincar\'e complexes with finite cyclic fundamental group of arbitrary order, and more generally with a finite fundamental group with $4$-periodic cohomology, are classified up to o.p.\ homotopy equivalence by their quadratic $2$-type.
Bauer~\cite{bauer} generalized this further to the case where only the $2$-Sylow subgroup of the fundamental group has $4$-periodic cohomology.

We consider the case where the $2$-Sylow subgroup is abelian with at most two generators. Note that abelian groups with at most two generators have $4$-periodic cohomology if and only if they are cyclic, so this improves previous results.

\begin{thm}
	\label{thm:main}
	Let $\pi$ be a finite group such that the $2$-Sylow subgroup
	is abelian with at most two generators.
	Then two oriented $4$-dimensional Poincar\'e complexes with fundamental group $\pi$
	are orientation-preserving homotopy equivalent if and only if their quadratic $2$-types are isomorphic.
\end{thm}

Our proof of \cref{thm:main} is based on the following result, which was proven by Hambleton-Kreck~\cite[Theorem~1.1~(i)]{HK}, with a hypothesis removed by Teichner~\cite{teichnerthesis} (see also~\cite[Corollary~1.5]{KT}).
Let $\Gamma$ be Whitehead's quadratic functor, whose definition we recall in \cref{sec:Gamma}.

\begin{thm}[{\cite{HK,teichnerthesis}}]
\label{thm:KT-cor-1.5}
  If $\Z\otimes_{\Z[\pi_1(X)]}\Gamma(\pi_2(X))$ is torsion free, then oriented $4$-dimensional Poincar\'e complexes with fundamental group $\pi$ are o.p.\ homotopy equivalent if and only if their quadratic $2$-types are isomorphic.
\end{thm}

By \cref{thm:KT-cor-1.5}, to prove \cref{thm:main} it suffices to show that $\Z\otimes_{\Z[\pi_1(X)]}\Gamma(\pi_2(X))$ is torsion free.  Here is an outline of our proof of this fact, which will occupy us for the remainder of the article, together with signposts as to where in the paper each step is carried out.
After recalling the definition and some properties of $\Gamma$, \cref{sec:Gamma} introduces tools for showing that groups of the form $\Z\otimes_{\Z\pi}\Gamma(L)$ are torsion free. \cref{sec:application-4-dim-cx} applies them to $L=\pi_2(X)$, where $X$ is a 4-dimensional Poincar\'{e} complex. This section recalls the short exact sequence for stable isomorphism classes of $\Z\pi$-modules $0 \to \ker d_2 \to \pi_2(X) \to \coker d^2 \to 0,$ where $d_2$ and $d^2$ come from a free $\Z\pi$-module resolution $(C_*,d_*)$ of $\Z$ and its dual respectively. At the start of \cref{sec:proof-main-thm}, we reduce the proof that $\Z\otimes_{\Z[\pi_1(X)]}\Gamma(\pi_2(X))$ is torsion free to showing that
$\mathcal{K} \coloneqq \Z\otimes_{\Z\pi} \Gamma(\ker d_2)$
and $\mathcal{CK} \coloneqq \Z\otimes_{\Z\pi} \Gamma(\coker d^2)$ are torsion free abelian groups, for $\pi$ a 2-Sylow subgroup of $\pi_1(X)$.
The remainder of \cref{sec:proof-main-thm}, which is the technical core of the paper, then proves that the groups $\mathcal{K}$ and $\mathcal{CK}$ are torsion free for $\pi$ any finite abelian group with at most two generators, completing the proof that $\Z\otimes_{\Z[\pi_1(X)]}\Gamma(\pi_2(X))$ is torsion free, and hence by \cref{thm:KT-cor-1.5} completing the proof of \cref{thm:main}.

\begin{rem}
Note that for nonorientable $4$-dimensional Poincar\'e complexes (or even manifolds) the analogous statement to \cref{thm:main} does not hold even for $\pi=\Z/2$, by \cite{KKR92}.  Nonetheless nonorientable closed 4-manifolds with fundamental group $\Z/2$ were classified up to homotopy equivalence by Hambleton, Kreck, and Teichner \cite{HKT-nonorientable}.
%

After this article appeared as a preprint, the first and third authors, together with Nicholson, proved the analogue of \cref{thm:main} for dihedral fundamental groups \cite{KNR}.
\end{rem}

We also note that the method of proof of \cref{thm:main} does not work for finite abelian groups whose $2$-Sylow subgroup requires more than two generators. As discussed in \cref{sebsection:groups-less-16}, for $\pi_1(X)\cong \Z/2 \times \Z/2 \times \Z/2$, the group $\Z \otimes_{\Z[\pi_1(X)]} \Gamma(\pi_2(X))$ is in general not torsion free. This leads us to pose the following question.

\begin{question}
Let $\pi$ be a finite abelian fundamental group whose $2$-Sylow subgroup requires more than two generators.
Are there oriented $4$-dimensional Poincar\'e complexes $($or even better, manifolds$)$ with fundamental group $\pi$ that are not homotopy equivalent but have isomorphic quadratic $2$-types?
\end{question}


\subsubsection*{Acknowledgments}
We are grateful to an anonymous referee for carefully reading the paper, and providing helpful comments that led to an improvement in the exposition.

The first author was funded by the Deutsche Forschungsgemeinschaft
under Germany's Excellence Strategy - GZ 2047/1, Projekt-ID 390685813.
The second author was partially supported by EPSRC New Investigator grant EP/T028335/2 and EPSRC New Horizons grant EP/V04821X/2.
The third author was supported by the Max Planck Institute for Mathematics in Bonn.

\section{Whitehead's \texorpdfstring{$\Gamma$}{Gamma} groups}
\label{sec:Gamma}

We recall the definition of Whitehead's $\Gamma$ functor from ~\cite{whitehead}, as well as a couple of key lemmas that we shall use in computations of $\Gamma$ groups.


\begin{defi}[$\Gamma$ groups]
	\label{def:gamma}
	Let $A$ be an abelian group. Then $\Gamma(A)$ is an abelian group with generators the elements of $A$. We write $a$ as $v(a)$ when we consider it as an element of $\Gamma(A)$. The group $\Gamma(A)$ has the following relations:
	\[\{v(-a)-v(a)\mid a\in A\}  \quad \text{ and }
	\]
	\[\{v(a+b+c)-v(b+c)-v(c+a)-v(a+b)+v(a)+v(b)+v(c)\mid a,b,c\in A\}.\]
\end{defi}
\noindent In particular, $v(0_A)=0_{\Gamma(A)}$.

\begin{rem}\label{rem:pi-2-free-abelian}
  For $X$ a $4$-dimensional Poincar\'{e} complex with finite fundamental group~$\pi$,
\[\pi_2(X) \cong H_2(\wt{X};\Z) \cong H^2(\wt{X};\Z) \cong \Hom(H_2(\wt{X};\Z),\Z),\]
which is in particular a finitely generated free abelian group.
\end{rem}

\begin{lem}[{\cite[page 62]{whitehead}}]
	\label{lem:gammafree}
	If $A$ is free abelian with basis $\mathfrak{B}$,
	then $\Gamma(A)$ is free abelian with basis
	\[
		\{ v(b), v(b+b')-v(b)-v(b') \mid b \neq b' \in \mathfrak{B} \}.
	\]
\end{lem}

In this case, we will usually consider $\Gamma(A)$
as the subgroup of symmetric elements of $A\otimes A$
given by sending $v(a)$ to $a\otimes a$.
Observe that $v(b+b')-v(b)-v(b')$ corresponds to
the symmetric tensor $b \otimes b' + b' \otimes b$.
For a $\Z\pi$-module $A$, the group $\pi$ acts on $\Gamma(A) \subseteq A \otimes A$ via
\[g\cdot \sum_i (a_i \otimes b_i) =\sum_i (g \cdot a_i) \otimes (g \cdot b_i).\]
To compute $\Gamma$ groups we will mostly rely on the following lemma.

\begin{lem}[{\cite[Lemma~4]{bauer}}]
	\label{lem:bauer}
Let $\pi$ be a group.	If $0\to A\to B\to C\to 0$ is a short exact sequence of $\Z\pi$-modules which are free as abelian groups, then there is a $\Z\pi$-module $D$, also free as an abelian group, and such that there are short exact sequences of $\Z\pi$-modules
	\[0\to \Gamma(A)\to \Gamma(B)\to D\to 0\]
	and
	\[0\to A\otimes_\Z C\to D\to \Gamma(C)\to 0.\]
\end{lem}

For the direct sum of free abelian groups $A, B$, the short exact sequences in the previous lemma split, so we have an isomorphism of $\Z\pi$-modules	
\[\Gamma(A \oplus B) \cong \Gamma(A) \oplus \Gamma(B) \oplus (A \otimes_{\Z} B).\]
For computational purposes we will need that the map $A\otimes_\Z C\to D$ is given as follows. Pick bases (as free abelian groups) $\{a_i\},\{c_j\}$ and $\{a_i,\wt c_j\}$ of $A$, $C$, and $B$ respectively, where $\wt c_j$ is a lift of $c_j$. Then $a_i\otimes c_j$ is sent to $[a_i\otimes \wt c_j +\wt c_j\otimes a_i] \in D \cong \Gamma(B)/\Gamma(A)$.

\subsection{Showing that \texorpdfstring{groups of the form $\Z\otimes_{\Z\pi}\Gamma(L)$}{tensored down Gamma groups} are torsion free}\label{sec:torsion-in-Gamma}
The strategy to show that groups of the form $\Z\otimes_{\Z\pi}\Gamma(L)$ are torsion free is to reduce to one of the cases in the next result, due to Hambleton and Kreck. See also their corrigendum \cite{HK-cor}, where a mistake in the original proof is fixed.

\begin{thm}[{\cite[Theorem~2.1]{HK}}]
	\label{thm:torsion-free-gamma}
	Let $\pi$ be a finite group. Let $L$ be one of the following:
	\begin{enumerate}[(i)]
		\item a finitely generated projective $\Z\pi$-module;
		\item the augmentation ideal $\ker(\varepsilon \colon \Z\pi \to \Z)$;
		\item  $\Z\pi/N$, where $N = \sum_{g \in \pi} g$ denotes the norm element.
	\end{enumerate}
	Then $\Z\otimes_{\Z\pi}\Gamma(L)$ is torsion free.
\end{thm}

We will also need the following proposition, which extends the third case of \cref{thm:torsion-free-gamma}.

\begin{prop}
	\label{prop:torsion-free-gamma2}
	Let $G$ and $H$ be finite groups and let $N_H=\sum_{h\in H}h$ be the norm element of $H$, considered in $\Z[G \times H]$.  Then $\Z\otimes_{\Z[G\times H]}\Gamma(\Z[G\times H]/N_H)$ is torsion free.
\end{prop}

For a $\Z H$-module $A$, let $A[G]:=\{\sum_{g\in G}a_gg\mid a_g\in A\}$
be the $\Z[G\times H]$-module
with the addition and the action given by
\[
	\big(\sum_{g\in G}a_gg\big) + \big(\sum_{g\in G}b_gg\big)
	=\sum_{g\in G}(a_g+b_g)g
\text{ and }
	(g',h)\big(\sum_{g\in G}a_gg\big)
	=\sum_{g\in G}(ha_g)(g'g)
\]
respectively. The proof of \cref{prop:torsion-free-gamma2} starts with the following lemma.

\begin{lem}
	\label{lem:torsion-free}
Let $G$ and $H$ be finite groups and write $\pi := G \times H$. Let $A$ be a $\Z H$-module which is finitely generated and free as an abelian group. Consider the subset $S'$ of $G$ of all elements with $g^2 \neq 1$. Then we have a free involution on $S'$ given by $g\mapsto g^{-1}$.  Let $S\subseteq G$ be a set of representatives for $S'/(\Z/2)$.
Then there is an isomorphism of $\Z\pi$-modules
\[
	\phi \colon \Gamma(A[G])
	\xrightarrow{\cong}
	\Gamma(A)[G] \oplus
	\bigoplus_{g\in S}(A\otimes_{\Z} A)[G] \oplus
	\bigoplus_{g\in G\setminus\{1\}, g^2=1}(A\otimes_{\Z} A)[G]/\textup{Flip}_g,
\]
where the submodule
$\textup{Flip}_g:=\langle (a\otimes b)-(b\otimes a)g\mid a,b\in A\rangle$.
\end{lem}

The special case $H=1$ and $A=\Z$ is \cite[Lemma~2.2]{HK} and the proof is similar.

\begin{proof}
The inverse of $\phi$ is given as follows. Elements $(a\otimes a')\gamma\in\Gamma(A)[G]$ are sent to $a\gamma\otimes a'\gamma$, elements $(a\otimes a')\gamma\in (A\otimes_\Z A)[G]$ in the factor indexed by $g\in S$ are sent to $a\gamma g\otimes a'\gamma+a'\gamma\otimes a\gamma g$, and for $g$ with $g^2=1$ elements $(a\otimes a')\gamma\in (A\otimes_\Z A)[G]/\langle (a\otimes b)-(b\otimes a)g\mid a,b\in A\rangle$ are again sent to $a\gamma g\otimes a'\gamma+a'\gamma\otimes a\gamma g$.

Since as an abelian group
\[\Gamma(A[G])\cong \bigoplus_{g\in \pi/H}\Gamma(A)\oplus \bigoplus_{\{g,g'\}\subseteq \pi/H,g\neq g'}A\otimes_{\Z} A,\]
it is easy to see that the above described map is surjective.
That it is also injective follows from a computation of the rank
of the involved modules considered as abelian groups, since both sides are free as abelian groups by \cref{lem:gammafree}.
\end{proof}

\begin{proof}[Proof of \cref{prop:torsion-free-gamma2}]
	Let $\pi:=G\times H$. We have $\Z[\pi]/N_H\cong (\Z H/N_H)[G]$. Thus by \cref{lem:torsion-free} the abelian group $\Z\otimes_{\Z\pi}\Gamma(\Z[\pi]/N_H)$ is a direct sum of groups of the form \[\Z\otimes_{\Z H} \Gamma(\Z H/N_H)\text{, } \Z\otimes_{\Z H}(\Z H/N_H\otimes_\Z \Z H/N_H) \text{ and } \Z\otimes_{\Z\langle H,g\rangle}(\Z H/N_H\otimes_\Z \Z H/N_H)\]
 where in the last case the element $g$ has order two and acts on $\Z H/N_H\otimes_\Z \Z H/N_H$ by flipping the two factors. The first of these is torsion free by \cref{thm:torsion-free-gamma}.

 As an abelian group,  $\Z H/N_H$ is free with $\Z$-basis given by $1\cdot h$ for $h\in H\setminus\{1\}$.
It is straightforward to check that the map
 \begin{align*}
   \Z \otimes_{\Z H} (\Z H/N_H \otimes_{\Z} \Z H/N_H) & \to \Z H/N_H ;\;\;
   1 \otimes ([\lambda] \otimes [\lambda'])  \mapsto [\overline{\lambda'} \lambda],
 \end{align*}
where the involution is determined by $\overline h=h^{-1}$,
is well-defined and an isomorphism of abelian groups.
It follows that $\Z\otimes_{\Z H}(\Z H/N_H\otimes_\Z \Z H/N_H)$ is a free abelian group.

Finally, the abelian group  $\Z\otimes_{\Z\langle H,g\rangle}(\Z H/N_H\otimes_\Z \Z H/N_H)$ can in the same way be identified with the orbits $(\Z H/N_H)/(\Z/2)$, where~$\Z/2$ acts by the involution.
This abelian group is also torsion free with $\Z$-basis
$(H\setminus \{1\})/(\Z/2)$,
where $\Z/2$ acts by inversion. This completes the proof of \cref{prop:torsion-free-gamma2}.
\end{proof}

\subsection{Connection with Tate homology}\label{sec:connection}

\begin{defi}[Tate homology]
	Given a finite group $\pi$ and a $\Z\pi$-module $A$, the Tate homology groups $\widehat H_n(\pi;A)$ are defined as follows. Let $\cdot N\colon A_\pi\to A^\pi$ denote multiplication with the norm element from the orbits $A_\pi:=\Z\otimes_{\Z\pi}A$ to the $\pi$-fixed points of $A$, that is $(\cdot N)(1 \otimes a) = \sum_{g \in \pi} g a \in A^\pi \subseteq A$.
Then
	\begin{align*}
	\widehat H_n(\pi;A) &:=H_n(\pi;A)\text{~for~}n\geq 1; \,\,
	\wh H_0(\pi;A) :=\ker(\cdot N); \\
	\wh H_{-1}(\pi;A)& :=\coker(\cdot N); \text{ and }
	\wh H_n(\pi;A) :=H^{-n-1}(\pi;A)\text{~for~}n\leq -2
	\end{align*}
\end{defi}

The following elementary observation will be central to our computations.

\begin{lem}\label{lemma:torsion-gamma-equals-tate}
	Let $\pi$ be a finite group.
	For a $\Z$-torsion free $\Z\pi$-module $A$, the torsion in $\Z\otimes_{\Z\pi} A$ is precisely the zeroth Tate homology $\wh H_0(\pi;A)$.
 \end{lem}

\begin{proof}
  Since the $\pi$-fixed points $A^\pi$ are again $\Z$-torsion free, and torsion maps to torsion under a homomorphism, the torsion in $\Z\otimes_{\Z\pi} A$ is contained in the kernel of multiplication by the norm element $N$. The composition of the multiplication with $N$ and the projection $A^\pi\to \Z\otimes_{\Z\pi} A$ is given by multiplication with the group order. Hence all elements in the kernel of multiplication with $N$ are annihilated by multiplication with $|\pi|$ and are thus torsion elements.  More precisely, suppose that $1 \otimes a \in \ker(\cdot N) \subseteq \Z \otimes_{\Z\pi} A$. Then in $\Z \otimes_{\Z\pi} A$ we have:
  \[|\pi|(1 \otimes a) = |\pi| \otimes a = 1 \otimes \sum_{g\in \pi} ga = 1 \otimes (\cdot N)(1\otimes a) = 1\otimes 0 =0. \]
So $1 \otimes a$ is a torsion element.
\end{proof}

\begin{lem}
	\label{lem:cor_of_bauer}
	For the sequences of coefficients in \cref{lem:bauer},
	$\Z\otimes_{\Z\pi}\Gamma(B)$ is torsion free if $\Z\otimes_{\Z\pi} \Gamma(A)$, $\Z\otimes_{\Z\pi}\Gamma(C)$ and $\Z\otimes_{\Z\pi}(A\otimes_\Z C)$ are torsion free.
\end{lem}
\begin{proof}
	Consider the long exact sequence of Tate homology groups~\cite[VI~(5.1)]{brown} corresponding to the sequence $0 \to \Gamma(A) \to \Gamma(B) \to D \to 0$,
	\[
		\cdots
		\to \wh H_1(\pi;D)
		\to \wh H_0(\pi; \Gamma(A))
		\to \wh H_0(\pi;\Gamma(B))
		\to \wh H_0(\pi;D)
		\to \wh H_{-1}(\pi;\Gamma(A))
		\to \cdots
	\]
	It follows that if $\Z\otimes_{\Z\pi}\Gamma(A)$ and $\Z\otimes_{\Z\pi}D$ are torsion free, i.e.\ $\wh H_0(\pi;\Gamma(A))=\wh H_0(\pi;D)=0$ by \cref{lemma:torsion-gamma-equals-tate}, then also $\Z\otimes_{\Z\pi} \Gamma(B)$ is torsion free. By the same argument, applied to the short exact sequence $0 \to A \otimes_{\Z} C \to D \to \Gamma(C)\to 0$, it follows from the assumptions that $\Z\otimes_{\Z\pi}\Gamma(C)$ and $\Z\otimes_{\Z\pi}(A\otimes_\Z C)$ are torsion free, that $\Z\otimes_{\Z\pi} D$ is torsion free.
\end{proof}


\section{\texorpdfstring{Application to $4$-dimensional Poincar\'{e} complexes}{Application to 4-dimensional Poincar\'{e} complexes}}\label{sec:application-4-dim-cx}

Let $X$ be a $4$-dimensional Poincar\'e complex
with finite fundamental group $\pi$.
Let $(C_*,d_*)$ be a free $\Z\pi$-resolution of $\Z$ with
$C_0=\Z\pi$, and with $C_1$ and $C_2$ finitely generated.
Then stably (possibly stabilising any of the three modules by a free module)
there is an extension
\[0\to \ker d_2\to \pi_2(X)\to \coker d^2\to 0,\]
as shown in~\cite[Proposition~2.4]{HK}.
Here the choice of resolution $(C_*,d_*)$ does not matter since for any
two choices of resolution the $\Z\pi$-modules
$\ker d_2$ and $\coker d^2$ are stably isomorphic;
see for example~\cite[Lemma~5.2]{KPT}, which relies on~\cite{lms197}.

\begin{lem}\label{lem:ker-and-coker-torsion-free}
  For every finite group $\pi$, and for every choice of free resolution $(C_*,d_*)$, $\ker d_2$ and $\coker d^2$ are torsion free abelian groups.
\end{lem}

\begin{proof}
Let  $(C_*,d_*)$ be a free $\Z\pi$-resolution of $\Z$ with
$C_0=\Z\pi$ and $C_1$ and $C_2$ finitely generated, as above.
Then $\ker d_2$ is a subgroup of a finitely generated free abelian group,
and so is torsion free.

Let $K$ be a finite $2$-dimensional CW complex with fundamental group $\pi$. Then the cellular chain complex $C_*(K;\Z\pi)$ is a start of a free resolution of $\Z$ as a $\Z\pi$-module. In particular $H^2(K;\Z\pi)$ is stably isomorphic to $\coker d^2$ and hence it suffices to show that $H^2(K;\Z\pi)$ is torsion free as an abelian group.
Since $\pi$ is finite, 	$H^2(K;\Z\pi)
	\cong H^2(\wt{K};\Z)$. By universal coefficients the latter group is isomorphic to $\Hom(H_2(\wt{K};\Z),\Z)$, which as required is torsion free.
\end{proof}

The next lemma combined with \cref{lem:ker-and-coker-torsion-free}
and \cref{lemma:torsion-gamma-equals-tate}
implies that the torsion subgroups of $\Z \otimes_{\Z\pi} \Gamma(\ker d_2)$
and $\Z \otimes_{\Z\pi} \Gamma(\coker d^2)$ only depend
on the stable isomorphism classes of $\ker d_2$ and $\coker d^2$ respectively,
and thus only depend on the group $\pi$.

\begin{lem}\label{lem:gamma-stable}
	Let $A$ be a $\Z\pi$-module that is free as an abelian group.
	Then $\wh H_0(\pi;\Gamma(A))$ only depends on the stable isomorphism type of $A$.
\end{lem}
\begin{proof}
	We have $\Gamma(A\oplus \Z\pi) \cong \Gamma(A) \oplus \Gamma(\Z\pi) \oplus (\Z\pi \otimes_\Z A).$
	By \cref{thm:torsion-free-gamma}, $\wh H_0(\pi;\Gamma(\Z\pi))=0$.
	Furthermore, $\wh H_0(\pi;\Z\pi\otimes_\Z A)$
	is the torsion in $\Z\pi\otimes_{\Z\pi} A\cong A$ which was torsion free by assumption. Hence $\wh H_0(\pi;\Gamma(A\oplus\Z\pi))\cong \wh H_0(\pi;\Gamma(A))$.
\end{proof}

We need the following elementary observation in the proof of the next proposition.

\begin{lem}[{\cite[III~(5.7)]{brown}}]\label{lem:free-module-check}
  Let $M$ be a left $\Z\pi$-module, consider $\Z\pi$ as a left rank one free module, and let $M \otimes_{\Z} \Z\pi$ be the $\Z\pi$-module where $\pi$ acts via the diagonal action. Let $(M \otimes_{\Z} \Z\pi)_r$ be the $\Z\pi$-module where $\pi$ acts just on $\Z\pi$ by left multiplication and $\pi$ acts trivially on $M$.  Then the map
  \begin{align*}
   \psi \colon M \otimes_{\Z} \Z\pi \to (M \otimes_{\Z} \Z\pi)_r ; \;\; m \otimes g  \mapsto g^{-1} m \otimes g
   \end{align*}
   is a $\Z\pi$-module isomorphism.
\end{lem}

\begin{proof}
 The map sending $m \otimes g \mapsto gm \otimes g$ gives an inverse. That both are $\Z\pi$-linear is straightforward to check.
\end{proof}

\begin{prop}\label{prop:always-torsion-free}
	The abelian group $\ker d_2\otimes_{\Z\pi} \coker d^2$
	is torsion free for every finite group $\pi$.
\end{prop}

\begin{proof}
This was noticed by Bauer \cite[page~5]{bauer}, but was not proven there.
We give the argument here. Consider the exact sequence
\[
	0 \to \Z
	\xrightarrow{N} C^0
	\xrightarrow{d^1} C^1
	\xrightarrow{d^2} C^2
	\to \coker d^2
	\to 0.
\]
By \cref{lem:free-module-check}, $\ker d_2\otimes_\Z C^i$ is $\Z\pi$-module isomorphic to $(\ker d_2\otimes_\Z C^i)_r$ the module with the same underlying abelian group but where the $\pi$ action is trivial on $\ker d_2$ and acts as usual on $C^i$.  Since $\ker d_2$ is free as a $\Z$-module, this is free as a $\Z\pi$-module.
Hence $\wh H_j(\pi;\ker d_2\otimes_\Z C^i)=0$ for all $i,j$ by \cite[VI~(5.3)]{brown}: if $A$ is a free $\Z\pi$-module, then $\wh H_j(\pi;A)=0$ for all $j\in\Z$.
%
%
Now, from the long exact sequence in Tate homology groups and dimension shifting~\cite[III.7]{brown} we get
	\[\wh H_0(\pi;\ker d_2\otimes_\Z \coker d^2)\cong \wh H_{-3}(\pi;\ker d_2\otimes_\Z \Z).\]
A similar argument shifting dimension upwards
	shows that $\wh H_{-3}(\pi;\ker d_2)\cong \wh H_0(\pi;\Z)=0.$
	Hence $\wh H_0(\pi;\ker d_2\otimes_\Z \coker d^2)=0$, which is equivalent to the statement of the proposition by \cref{lemma:torsion-gamma-equals-tate}, since $\Z \otimes_{\Z\pi} (\ker d_2\otimes_{\Z} \coker d^2) \cong \ker d_2\otimes_{\Z\pi} \coker d^2$.
\end{proof}


\begin{cor}
	\label{cor:torsionfree}
	Assume that $\Z\otimes_{\Z\pi} \Gamma(\ker d_2)$ and $\Z\otimes_{\Z\pi} \Gamma(\coker d^2)$ are torsion free. Then $\Z\otimes_{\Z\pi} \Gamma(\pi_2(X))$ is torsion free.
\end{cor}

\begin{proof}
	The torsion in $\Z\otimes_{\Z\pi} \Gamma(\pi_2(X))$ equals
	$\wh H_0(\pi;\Gamma(\pi_2(X))$ by \cref{lem:ker-and-coker-torsion-free} and \cref{lemma:torsion-gamma-equals-tate}.
	By \cref{lem:ker-and-coker-torsion-free} again and \cref{lem:gamma-stable}, this torsion only depends on the stable isomorphism class
	of $\pi_2(X)$, and therefore we may consider the stable exact sequence
	$0\to \ker d_2\to \pi_2(X)\to \coker d^2\to 0$.
	By \cref{lem:cor_of_bauer}, $\Z\otimes_{\Z\pi} \Gamma(\pi_2(X))$ is torsion free
	if $\Z\otimes_{\Z\pi} \Gamma(\ker d_2)$, $\Z\otimes_{\Z\pi}\Gamma(\coker d^2)$
	and $\Z\otimes_{\Z\pi}(\ker d_2\otimes_\Z \coker d^2)$ are torsion free.
	The first two hold by assumption, the latter holds by
	\cref{prop:always-torsion-free}.
\end{proof}

As stated in the introduction, by \cite[Corollary~1.5]{KT}, to show \cref{thm:main} it suffices to show that $\Z\otimes_{\Z[\pi_1(X)]}\Gamma(\pi_2(X))$ is torsion free.
We will therefore want to show that $\Z\otimes_{\Z\pi} \Gamma(\ker d_2)$ and $\Z\otimes_{\Z\pi} \Gamma(\coker d^2)$ are torsion free in the cases we consider.  We will use the following lemma several times.

\begin{lem}
	\label{lem:divisible}
	Let $\pi$ be a finite group. Let $0\to A\to B\to C\to 0$ be a short exact sequence of $\Z\pi$-modules such that $B$ is torsion free as a $\Z$-module and such that $\Z\otimes_{\Z\pi}A$ is torsion free. In this case the sequence
	\[
		0 \to \Z\otimes_{\Z\pi} A
		\to \Z\otimes_{\Z\pi} B
		\to \Z\otimes_{\Z\pi} C
		\to 0
	\]
	is exact. Assume that for every $($nontrivial$)$ torsion element $t\in \Z\otimes_{\Z\pi} C$, of order $k$ say, there is a preimage $b\in \Z \otimes_{\Z\pi} B$ of $t$ such that the preimage $a \in \Z \otimes_{\Z\pi} A$ of $kb \in \Z \otimes_{\Z\pi} B$ is not divisible by $k$.
	Then $\Z\otimes_{\Z\pi} B$ is a torsion free $\Z$-module.
\end{lem}

\begin{proof}
	We have an exact sequence
	\[H_1(\pi;C)\to \Z\otimes_{\Z\pi} A\to \Z\otimes_{\Z\pi} B\to \Z\otimes_{\Z\pi} C\to 0.\]
	As $\pi$ is finite, $H_1(\pi;C)$ is a torsion group~\cite[III~(10.2)]{brown}. Since $\Z\otimes_{\Z\pi} A$ is torsion free, the map $H_1(\pi;C)\to \Z\otimes_{\Z\pi} A$ is trivial and thus the sequence
	\[0\to \Z\otimes_{\Z\pi} A\to \Z\otimes_{\Z\pi} B\to \Z\otimes_{\Z\pi} C\to 0\]
	is exact.
	
	Since $\Z\otimes_{\Z\pi} A$ is torsion free, every nontrivial torsion element in $\Z\otimes_{\Z\pi} B$ has to map to a nontrivial torsion element in $\Z\otimes_{\Z\pi} C$. Now suppose for a contradiction that $\Z \otimes_{\Z\pi} B$ has torsion.  Let $t'\in \Z\otimes_{\Z\pi} B$ be a nontrivial torsion element mapping to $t\in  \Z\otimes_{\Z\pi} C$ of order $k$. Let $b\in \Z \otimes_{\Z\pi} B$ be a preimage as in the assumption, and let $a \in \Z\otimes_{\Z\pi} A$ be the preimage of $kb$, where by assumption~$a$ is not divisible by~$k$. Then $b-t'$ maps to $0 \in \Z \otimes_{\Z\pi} C$. Let $a' \in \Z\otimes_{\Z\pi} A$ be the preimage of $b-t'$. Furthermore,~$kt'$ maps to $kt=0 \in \Z\otimes_{\Z\pi} C$ and is torsion, thus $kt'$ is trivial in $\Z \otimes_{\Z\pi} B$. Hence~$ka'$ maps to $k(b-t') = kb$. Thus~$a=ka'$. But~$ka'$ is divisible by~$k$, contradicting the assumption on~$a=ka'$. It follows that $\Z\otimes_{\Z\pi} B$ is a torsion free $\Z$-module, as desired.
\end{proof}

\begin{rem}
	\label{rem:torsion}
	Suppose the torsion in $\Z \otimes_{\Z \pi} C$ is
	generated by a single element $x$ of order $n$.
	The following observation allows us to
	simplify what we have to check to apply the second
	half of \cref{lem:divisible}.
	The order of the multiple $t = \ell x$
	is $\frac{n}{\gcd(n, \ell)}$.
	Let $b \in \Z \otimes_{\Z \pi} B$ be a preimage of
	$x$. Then $\ell b$ is a preimage of $t=\ell x$. If
	the preimage $a'$ of $\frac{n}{\gcd(n, \ell)} \ell b$
	in $\Z \otimes_{\Z \pi} A$ is divisible
	by $\frac{n}{\gcd(n, \ell)}$,
	then the preimage $a$ of $n b$ would also be divisible by $\frac{n}{\gcd(n, \ell)}$
	(since $\frac{n}{\gcd(n,\ell)}$ is coprime to $\frac{\ell}{\gcd (n,\ell)}$).
	Hence to show that $\Z \otimes_{\Z \pi} B$
	is torsion free in the case
	that the torsion in
	$\Z \otimes_{\Z \pi} C$ is cyclic, it suffices to
	check that the preimage $a$ is not divisible by any
	divisor of $n$.	
\end{rem}

\section{Proof of the main theorem}\label{sec:proof-main-thm}

\begin{proof}[Proof of \cref{thm:main}]
Recall that by \cite[Corollary~1.5]{KT}, to show \cref{thm:main} it suffices to show that $\Z\otimes_{\Z[\pi_1(X)]}\Gamma(\pi_2(X))$ is torsion free.

\cref{cor:torsionfree} implies that it is enough to show that $\Z\otimes_{\Z\pi} \Gamma(\ker d_2)$ and $\Z\otimes_{\Z\pi} \Gamma(\coker d^2)$ are torsion free, where $(C_*,d_*)$ is a free $\Z\pi$-resolution of $\Z$ with
$C_0=\Z\pi$, and $C_1$ and $C_2$ finitely generated.  By \cref{lemma:torsion-gamma-equals-tate}, we therefore have to show that
\[
	\wh H_0(\pi;\Gamma(\ker d_2)) = 0
	\text{ and }
	\wh H_0(\pi;\Gamma(\coker d^2)) = 0.
\]
As observed by Bauer \cite{bauer},
it suffices to show this for a
2-Sylow subgroup $G$ of $\pi$.
For every non-cyclic finite abelian group $\pi$ with two generators, we shall show in \cref{prop:ker} and \cref{prop:coker} below that $\wh H_0(\pi;\Gamma(\ker d_2))$ and $\wh H_0(\pi;\Gamma(\coker d^2))$ vanish.  In particular this holds for $\pi$ the two generator abelian 2-Sylow subgroup of some larger group.

We also need to check 1-generator finite abelian groups, i.e.\ finite cyclic groups. This was done by Hambleton and Kreck \cite{HK} but we repeat the argument here.
Let $C$ be the cyclic group of order~$k$ generated by $t$ and as ever let $N = \sum_{i=0}^{k-1} t^i$ be the norm element.
There is a free resolution whose first few terms are
\[
	\cdots \to
	\Z C \xrightarrow{1-t}
	\Z C \xrightarrow{\cdot N}
	\Z C \xrightarrow{1-t} \Z C
	\to \Z
	\to 0.
\]
Therefore $\ker d_2 \cong \langle 1-t\rangle \cong \Z C/N$
and $\coker d^2 \cong \Z C/N$. By \cref{thm:torsion-free-gamma}~(iii),
we have that $\Z \otimes_{\Z C}\Gamma(\Z C/N)$ is torsion free as required.
Modulo \cref{prop:ker} and \cref{prop:coker} below, this completes the proof of \cref{thm:main}.
\end{proof}

\subsection{The computation \texorpdfstring{of $\wh H_0(\pi;\Gamma(\ker d_2))$}{for the kernel}}\label{sec:computation-kernel}

In this section we show the following.
\begin{prop}
	\label{prop:ker}
For every finite abelian group with two generators $\pi$, $\wh H_0(\pi;\Gamma(\ker d_2))=0$.
\end{prop}

As explained above, we may assume that $|\pi|$ is a power of 2 for the purpose of the proof of \cref{thm:main}, but we do not need this assumption for \cref{prop:ker}.

\begin{proof}
  For the group $\pi=\langle a,b\mid a^n,b^m,[a,b]\rangle$ let $N_a:=\sum_{i=0}^{n-1}a^i$
and $N_b:=\sum_{i=0}^{m-1}b^i$.
Let $C_2\xrightarrow{d_2}C_1\xrightarrow{d_1}C_0$
be the chain complex  of $\Z\pi$-modules corresponding to the presentation
$\langle a,b\mid a^n,b^m,[a,b]\rangle$.
Extend this to the standard free resolution
of $\Z$ as a $\Z\pi$-module:
\[
\xymatrix@C-0.35cm{
	C_{4} \ar[d]^-{d_{4}} &
	\Z\pi \ar[dr]^-{N_{a}} & \oplus &
	\Z\pi \ar[dl]^-{b-1} \ar[dr]^-{1-a} & \oplus &
	\Z\pi \ar[dl]^-{N_{b}} \ar[dr]^-{N_{a}} & \oplus &
	\Z\pi \ar[dl]^-{b-1} \ar[dr]^-{1-a} & \oplus &
	\Z\pi \ar[dl]^-{N_{b}} \\
	C_{3} \ar[d]^-{d_{3}} &
	& \Z\pi \ar[dr]^-{1-a} & \oplus &
	\Z\pi \ar[dl]^-{1-b} \ar[dr]^{N_{a}} & \oplus &
	\Z\pi \ar[dl]^-{-N_{b}} \ar[dr]^-{1-a} & \oplus &
	\Z\pi \ar[dl]^-{1-b}
	\\
	C_{2} \ar[d]^-{d_{2}} &
	& & \Z\pi \ar[dr]^-{N_{a}} & \oplus &
	\Z\pi \ar[dl]^-{b-1} \ar[dr]^-{1-a} & \oplus &
	\Z\pi \ar[dl]^-{N_{b}}
	\\
	C_{1} \ar[d]^-{d_{1}} &
	& & & \Z\pi \ar[dr]^-{1-a} & \oplus &
	\Z\pi \ar[dl]^-{1-b}
	\\
	C_{0} &
	& & & & \Z\pi
  }
\]
By exactness, $\ker d_2\cong \operatorname{im} d_3\cong C_3/\ker d_3 \cong \coker d_4$. From this it follows that
\[
	\ker d_2\cong (\Z\pi)^4/\langle (N_a,0,0,0),(b-1,1-a,0,0),(0,N_b,N_a,0), (0,0,b-1,1-a),(0,0,0,N_b)\rangle.
\]
Define
\[M_1:=(\Z\pi)^2/\langle (N_a,0),(0,N_b)\rangle
\text{ and }
M_2:=(\Z\pi)^2/\langle (1-a,0),(N_b,N_a),(0,b-1) \rangle.\]
The inclusion of the outer two summands of $(\Z\pi)^4$ induces a short exact sequence
\begin{equation}\label{eq:K}0\to M_1 \to \ker d_2\to M_2 \to 0.\end{equation}
For $M_2$ we have a short exact sequence
\begin{equation}
\label{eq:M2}
0\to \Z\pi/\langle 1-a\rangle\to M_2\to \Z\pi/\langle N_a,b-1\rangle\to 0.
\end{equation}
In both these sequences, a quick check is required that the map on the left is indeed injective.

\begin{claim}
	\label{claim:prop_kernel}
	We have that $\Z \otimes_{\Z\pi}\Gamma(M_2)$ is torsion free,
	that is $\wh{H}_0(\pi;\Gamma(M_2))=0$.
\end{claim}
The proof of this claim takes more than a page.
To prove the claim, we want to apply \cref{lem:divisible} to the short exact sequence
\begin{equation}\label{eq:M2-ses-1}
0 \to \Gamma(\Z\pi/\langle 1-a\rangle) \to \Gamma(M_2) \to D \to 0\end{equation}
coming from~\cref{lem:bauer}.
We will use the other sequence
\begin{equation}\label{eq:M2-ses-2}
0 \to \Z\pi/\langle 1-a\rangle \otimes_{\Z} \Z\pi/\langle N_a,b-1\rangle \to D \to \Gamma(\Z\pi/\langle N_a,b-1\rangle) \to 0,\end{equation}
from \cref{lem:bauer} as well, to compute the torsion in $\wh{H}_0(\pi;D)$.

We check that \eqref{eq:M2-ses-1} satisfies the hypotheses of \cref{lem:divisible} that
$\Gamma(M_2)$ and $\Z \otimes_{\Z\pi} \Gamma(\Z\pi/\langle 1-a\rangle)$ are torsion free as $\Z$-modules.
Let $C_a,C_b$ denote the cyclic subgroups of $\pi$ generated by $a$ and $b$ respectively, so $\pi \cong C_a \times C_b$.
Note that as abelian groups $\Z\pi/\langle 1-a\rangle \cong \Z C_b$ and $\Z\pi/\langle N_a,b-1\rangle \cong \Z C_a / N_a$.
In particular both are torsion free abelian groups.
It follows from the extension \eqref{eq:M2}
that $M_2$ is also torsion free,
then \cref{lem:gammafree} implies that $\Gamma(M_2)$ is torsion free.
Furthermore, we have
\[\wh H_0(\pi;\Gamma(\Z\pi/\langle 1-a\rangle))\cong \wh H_0(\pi;\Gamma(\Z C_b)) \cong \wh H_0(C_b;\Gamma(\Z C_b)).\]
The second isomorphism follows by observing that
$\wh{H}_0(\pi;\Gamma(\Z C_b))$ is the torsion in $\Z \otimes_{\Z\pi} \Gamma(\Z C_b)$,
while $\wh{H}_0(C_b;\Gamma(\Z C_b))$ is the torsion in $\Z \otimes_{\Z C_b} \Gamma(\Z C_b)$.
But the action of $\Z\pi$ on $\Gamma(\Z C_b)$ factors through the homomorphism $\Z\pi \to \Z\pi/\langle 1-a\rangle \xrightarrow{\cong} \Z C_b$, so tensoring over $\Z\pi$ and tensoring over $\Z C_b$ yield isomorphic groups.  Now, $\wh H_0(C_b;\Gamma(\Z C_b))=0$ by \cref{lemma:torsion-gamma-equals-tate} and \cref{thm:torsion-free-gamma}~(i), so $\Z \otimes_{\Z\pi} \Gamma(\Z\pi/\langle 1-a\rangle)$ is torsion free.  This completes the proof that \eqref{eq:M2-ses-1} satisfies the hypotheses from the first sentence of \cref{lem:divisible} that $\Gamma(M_2)$ and $\Z \otimes_{\Z\pi} \Gamma(\Z\pi/\langle 1-a\rangle)$ are torsion free as $\Z$-modules.

We continue proving the claim that $\Z \otimes_{\Z\pi}\Gamma(M_2)$ is torsion free.
The first part of \cref{lem:divisible} now gives us a short exact sequence
\begin{equation}\label{eq:M2-computation-Z-seq-1}
 0 \to \Z \otimes_{\Z\pi} \Gamma(\Z\pi/\langle 1-a\rangle) \to \Z \otimes_{\Z\pi} \Gamma(M_2) \to \Z \otimes_{\Z\pi} D \to 0.
 \end{equation}
To apply this sequence, we need to understand the torsion in $\Z \otimes_{\Z\pi} D$, which we do next.

By tensoring over $\Z\pi$ with $\Z$, the sequence \eqref{eq:M2-ses-2} gives rise to a long exact sequence ending in:
\[
	\cdots \to \Z\otimes_{\Z\pi} (\Z\pi/\langle 1-a\rangle\otimes_{\Z}\Z\pi/\langle N_a,b-1\rangle)
	\to \Z\otimes_{\Z\pi} D
	\to \Z\otimes_{\Z\pi} \Gamma(\Z\pi/\langle N_a,b-1\rangle)
	\to 0.
\]
We noted above that $\Z\pi/\langle N_a,b-1\rangle \cong \Z C_a / N_a$. Since the action of $\Z\pi$ on this module factors through $\Z\pi \to \Z\pi/\langle b-1 \rangle \xrightarrow{\cong} \Z C_a$, the same holds for the $\Z\pi$ action on $\Gamma(\Z C_a / N_a)$, and so $\Z \otimes_{\Z\pi} \Gamma(\Z C_a / N_a) \cong  \Z \otimes_{\Z C_a} \Gamma(\Z C_a / N_a)$. This latter group is torsion free by \cref{thm:torsion-free-gamma}~(iii).
Therefore the torsion in $\Z\otimes_{\Z\pi} D$ comes from
\begin{align*}
\Z\otimes_{\Z\pi} (\Z\pi/\langle 1-a\rangle\otimes_{\Z}\Z\pi/\langle N_a,b-1\rangle) &\cong \Z\pi/\langle 1-a\rangle\otimes_{\Z\pi} \Z\pi/\langle N_a,b-1\rangle \\  &\cong \Z\pi/\langle 1-a,N_a,b-1\rangle\cong \Z/n.
\end{align*}

The image of the generator $1\otimes 1$ of $\Z\pi/\langle 1-a\rangle\otimes_{\Z\pi}\Z\pi/\langle N_a,b-1\rangle$ in $\Z\otimes_{\Z\pi} D$ has as a preimage in $\Z\otimes_{\Z\pi} \Gamma(M_2)$
the element $[(1,0)\otimes (0,1)+(0,1)\otimes(1,0)]$.  Here we represent elements of $\Z\otimes_{\Z\pi} \Gamma(M_2)$ as symmetric tensors in $(\Z\pi)^2 \otimes_{\Z} (\Z\pi)^2$: each element of $(\Z\pi)^2$ determines an element of $M_2$, so we obtain an element of $\Gamma(M_2) \subseteq M_2\otimes_{\Z} M_2$. The square brackets around $[(1,0)\otimes (0,1)+(0,1)\otimes(1,0)]$ indicates taking orbits under the $\pi$ action i.e.\ the element $1 \otimes ((1,0)\otimes (0,1)+(0,1)\otimes(1,0)) \in \Z \otimes_{\Z\pi} \Gamma(M_2)$.
From now on we will use this square bracket notation to denote the classes of elements in various quotients of tensor products of free $\Z\pi$-modules.

Since the generator of $\Z\pi/\langle 1-a\rangle\otimes_{\Z\pi}\Z\pi/\langle N_a,b-1\rangle$ has order $n$, by exactness of \eqref{eq:M2-computation-Z-seq-1}, we have that
 \[[n((1,0)\otimes (0,1)+(0,1)\otimes(1,0))]\in \Z\otimes_{\Z\pi}\Gamma(M_2)\] lies in the image of $\Z\otimes_{\Z\pi}\Gamma(\Z\pi/\langle 1-a\rangle)$, in the sequence~\eqref{eq:M2-computation-Z-seq-1}.
To apply \cref{lem:divisible} together with \cref{rem:torsion},
we have to show that $[n((1,0)\otimes (0,1)+(0,1)\otimes(1,0))]$ has preimage in
$\Z\otimes_{\Z\pi}\Gamma(\Z\pi/\langle 1-a\rangle)$
that is not divisible by $n$ nor any of its factors.
Let $\overline{\cdot}$ be the usual involution sending $g \mapsto g^{-1}$.
We find the preimage of $[n((1,0)\otimes (0,1)+(0,1)\otimes(1,0))]$:
\begin{align*}
[n((1,0)\otimes (0,1)+(0,1)\otimes(1,0))]&=[(N_a,0)\otimes (0,1)+(0,1)\otimes(N_a,0)] \\
= \sum_{i=0}^{n-1} \big[(a^i,0) \otimes (0,1) + (0,1) \otimes (a^i,0)\big]  &= \sum_{i=0}^{n-1} \big[(1,0) \otimes (0,a^{-i}) + (0,a^{-i}) \otimes (1,0)\big] \\
=[(1,0)\otimes(0,\overline{N}_a)+(0,\overline{N}_a)\otimes(1,0)] & =[(1,0)\otimes(0,N_a)+(0,N_a)\otimes(1,0)] \\ &=[-(1,0)\otimes (N_b,0)-(N_b,0)\otimes(1,0)].
\end{align*}
Here the first equation uses that $a$ acts trivially on the first factor,
so multiplication by $n$ and by $N_a$ are equivalent.
The second and fourth equations use the definition of $N_a$.
The third equation uses that we have tensored with $\Z$ over $\Z\pi$.
In general, if $\Z \pi$ acts diagonally on a tensor
product $L \otimes_{\Z} L$ of $\Z \pi$-modules,
in the tensored down module
$\Z \otimes_{\Z \pi} (L \otimes L)$ the
relation $[a \otimes (\lambda b)] = [(\overline{\lambda} a) \otimes b]$
holds, where $\lambda \in \Z \pi$ and $a,b \in L$.
The fifth equation uses that $\overline{N}_a=N_a$.
The sixth and final equation uses the second relation of $M_2$.

The preimage element $1\otimes N_b+N_b\otimes 1$ is not divisible in $\Z\otimes_{\Z\pi}\Gamma(\Z\pi/\langle 1-a\rangle)\cong \Z\otimes_{\Z C_b}\Gamma(\Z C_b)$.
This follows from the concrete description of $\Gamma(\Z C_b)$ as a $\Z C_b$-module from \cite[Lemma~2.2]{HK}; see also \cref{lem:torsion-free}.
Now by \cref{lem:divisible} and \cref{rem:torsion}, $\Z \otimes_{\Z\pi}\Gamma(M_2)$ is torsion free, so $\wh H_0(\pi;\Gamma(M_2))$ is trivial.
This completes the proof of the claim.
\newline

Now we show that $\Z\otimes_{\Z\pi} \Gamma(\ker d_2)$
is torsion free, as desired for the proposition.
We want to apply the short exact sequence~\eqref{eq:K}.
By \cref{lem:bauer} we have a $\Z\pi$-module $E$ and short exact sequences
\begin{equation}\label{eq:M1-ses-1}
  0 \to \Gamma(M_1) \to \Gamma(\ker d_2) \to E \to 0
\end{equation}
and
\begin{equation}\label{eq:M1-ses-2}
0 \to M_1 \otimes_{\Z} M_2 \to E \to \Gamma(M_2) \to 0.
\end{equation}

\begin{claim}
  The short exact sequence~\eqref{eq:M1-ses-1} satisfies the hypotheses from the first sentence of \cref{lem:divisible}. That is, $\Gamma(\ker d_2)$ and $\Z \otimes_{\Z\pi} \Gamma(M_1)$ are torsion free as $\Z$-modules.
\end{claim}

We have $M_1\cong \Z\pi/ N_a \oplus \Z\pi/ N_b$ and hence $M_1$ is torsion free. We saw above, just below \eqref{eq:M2-ses-2}, that $M_2$ is torsion free.
Then the short exact sequence $0 \to M_1 \to \ker d_2 \to M_2 \to 0$, or alternatively \cref{lem:ker-and-coker-torsion-free}, implies that $\ker d_2$ is torsion free. Therefore $\Gamma(\ker d_2)$ is torsion free by \cref{lem:gammafree}.
Next, $M_1\cong \Z\pi/ N_a \oplus \Z\pi/ N_b$ implies that
\[\Gamma(M_1)\cong \Gamma(\Z\pi/ N_a)\oplus \Gamma(\Z\pi/ N_b)\oplus (\Z\pi/ N_a \otimes_\Z \Z\pi/ N_b).\]
The groups  $\Gamma(\Z\pi/ N_a)$ and $\Gamma(\Z\pi/ N_b)$ are still torsion free after applying $\Z\otimes_{\Z\pi}$ by \cref{prop:torsion-free-gamma2}.  The group
\[
	\Z \otimes_{\Z\pi} (\Z\pi/N_a \otimes_\Z \Z\pi/ N_b)
	\cong \Z\pi/N_a \otimes_{\Z\pi} \Z\pi/N_b
	\cong \Z\pi/\langle N_a,N_b\rangle
\]
is torsion free with $\Z$-basis given by $\{a^ib^j\mid 1\leq i\leq n-1, 1\leq j\leq m-1\}$. Therefore $\Z \otimes_{\Z\pi} \Gamma(M_1)$ is torsion free.  This completes the proof of the claim that \eqref{eq:M1-ses-1} satisfies the hypotheses of the first sentence of \cref{lem:divisible}.
\newline

Next we compute the torsion in $\Z \otimes_{\Z\pi} E$, i.e.\ $\wh H_0(\pi;E)$. The short exact sequence~\eqref{eq:M1-ses-2} gives rise to a long exact sequence:
\[\cdots\to \wh H_0(\pi;M_1\otimes_\Z M_2)\to \wh H_0(\pi;E)\to \wh H_0(\pi;\Gamma(M_2))\to\cdots\]
By the claim above, $\wh H_0(\pi;\Gamma(M_2))=0$.
We have
\[
	\Z \otimes_{\Z\pi} (M_1\otimes_{\Z} M_2)\cong M_1\otimes_{\Z\pi} M_2 \cong
	(\Z\pi/N_a \otimes_{\Z\pi} M_2)
	\oplus
	(\Z\pi/N_b \otimes_{\Z\pi} M_2).
\]
Using the description of $M_2$ in terms of generators and relations, \[M_2:=(\Z\pi)^2/\langle (1-a,0),(N_b,N_a),(0,b-1) \rangle,\] it follows that
\begin{align*}
\Z\pi/N_a\otimes_{\Z\pi}M_2
&\cong (\Z\pi)^2/\langle(1-a,0),(N_b,0),(0,b-1),(N_a,0),(0,N_a)\rangle\\
&\cong (\Z/n)C_b/N_b\oplus \Z C_a/N_a.
\end{align*}
Similarly, \[\Z\pi/N_b\otimes_{\Z\pi}M_2\cong \Z C_b/N_b\oplus (\Z/m) C_a/N_a.\]
It follows that the torsion in $\Z\otimes_{\Z\pi}E$ has preimages in $\Z\otimes_{\Z\pi}\Gamma(\ker d_2)$ of the form
\begin{align*}
[(1,0,0,0)\otimes (0,\lambda,0,0) + (0,\lambda,0,0)\otimes(1,0,0,0)  +(0,0,0,1)\otimes (0,0,\lambda',0)+(0,0,\lambda',0)\otimes(0,0,0,1)]\end{align*} with $\lambda\in \Z C_b$ and $\lambda'\in \Z C_a$.
Here we represent elements of $\Z\otimes_{\Z\pi} \Gamma(\ker d_2)$ as symmetric tensors in $(\Z\pi)^4 \otimes_{\Z} (\Z\pi)^4$, which determine elements in $\Gamma(\ker d_2) \subseteq \ker d_2 \otimes_{\Z} \ker d_2$. As above square brackets indicate taking orbits under the $\pi$ action, that is the corresponding element of $\Z \otimes_{\Z\pi} \Gamma(\ker d_2)$.

\begin{claim}
  The torsion in $\Z\otimes_{\Z\pi}E$ is generated by the elements
  \[[(1,0,0,0)\otimes (0,1,0,0) + (0,1,0,0)\otimes(1,0,0,0)] \text{ and } [(0,0,0,1)\otimes (0,0,1,0)+(0,0,1,0)\otimes(0,0,0,1)],\]
  that is, elements corresponding to the first two summands in the sum above with $\lambda=1$, and the final two summands with $\lambda'=1$.
\end{claim}

Now we prove the claim. Using the relation $(b-1,1-a,0,0)$ in $\ker d_2$,
in $\ker d_2\otimes_{\Z}\ker d_2$ we have \[(1-b,0,0,0)\otimes(1-b,0,0,0)=(0,1-a,0,0)\otimes(0,1-a,0,0).\] Hence the same relation holds in $\Gamma(\ker d_2)$. We have
\begin{align*}
  (0,1-a,0,0)\otimes(0,1-a,0,0)  = &(0,1,0,0)\otimes(0,1,0,0)-(0,a,0,0)\otimes(0,1,0,0) \\ &- (0,1,0,0)\otimes(0,a,0,0)+(0,a,0,0)\otimes(0,a,0,0).
\end{align*}
In $\Z\otimes_{\Z\pi}\Gamma(\ker d_2)$ the last element represents
the same element as $(0,1,0,0)\otimes(0,1,0,0)$
since we can act diagonally with $a^{-1}$.
Hence in $\Z\otimes_{\Z\pi}\Gamma(\ker d_2)$ we have
\[[(0,1-a,0,0)\otimes(0,1-a,0,0)]=[(0,1-a,0,0)\otimes(0,1,0,0)+(0,1,0,0)\otimes(0,1-a,0,0)].\]
Using the same relation as before, we have
\begin{align*}
	&[(0,1-a,0,0)\otimes(0,1,0,0)+(0,1,0,0)\otimes(0,1-a,0,0)] \\
	=&[(1-b,0,0,0)\otimes(0,1,0,0)+(0,1,0,0)\otimes(1-b,0,0,0)]\\
	=&[(1,0,0,0)\otimes(0,1-b^{-1},0,0)+(0,1-b^{-1},0,0)\otimes(1,0,0,0)].
\end{align*}
Using that $(1-b,0,0,0)\otimes(1-b,0,0,0)\in \Gamma(M_1)$, in the quotient $\Z\otimes_{\Z\pi} E$ this element is trivial. Hence in $\Z\otimes_{\Z\pi} E$ we have
\begin{align*}
	[(1,0,0,0)\otimes(0,\lambda,0,0)+(0,\lambda,0,0)\otimes(1,0,0,0)] =[(1,0,0,0)\otimes(0,|\lambda|,0,0)+(0,|\lambda|,0,0)\otimes(1,0,0,0)],
\end{align*}
where $|\lambda|=\sum_{i=0}^{m-1}\lambda_{b^i}\in\Z$ for $\lambda=\sum_{i=0}^{m-1}\lambda_{b^i}b^i$.
An analogous argument shows that
\begin{align*}
	&[(0,0,0,1)\otimes (0,0,\lambda',0)+(0,0,\lambda',0)\otimes(0,0,0,1)]\\
	=&[(0,0,0,1)\otimes (0,0,|\lambda'|,0)+(0,0,|\lambda'|,0)\otimes(0,0,0,1)],
\end{align*}
in the quotient $\Z\otimes_{\Z\pi} E$, where $|\lambda'| = \sum_{i=0}^{n-1}\lambda'_{a^i}\in \Z$ for $\lambda'=\sum_{i=0}^{m-1}\lambda'_{a^i}a^i$.
This shows that the torsion in $\Z\otimes_{\Z\pi} E$ is generated by $[(1,0,0,0)\otimes (0,1,0,0)+(0,1,0,0)\otimes(1,0,0,0)]$ and $[(0,0,0,1)\otimes (0,0,1,0)+(0,0,1,0)\otimes(0,0,0,1)]$, which completes the proof of the claim.
\newline

We continue proving that $\Z\otimes_{\Z\pi} \Gamma(\ker d_2)$
is torsion free, by applying \cref{lem:divisible}.
Let $x_a:=\sum_{i=1}^{n-1}(n-i)a^i$. A short computation shows that $n-N_a = x_a(a^{-1}-1)$.
Then using the relations $(N_a,0,0,0)$ and $(b-1,1-a,0,0)$ in $\ker d_2$, we have:
\begin{align*}
&[n((1,0,0,0)\otimes (0,1,0,0)+(0,1,0,0)\otimes(1,0,0,0))]\\
=&[(n-N_a,0,0,0)\otimes (0,1,0,0)+ (0,1,0,0)\otimes(n-N_a,0,0,0)]\\
=&[(-x_a(1-a^{-1}),0,0,0)\otimes (0,1,0,0)+(0,1,0,0)\otimes(-x_a(1-a^{-1}),0,0,0)]\\
=&[(-x_a,0,0,0)\otimes (0,1-a,0,0)+(0,1-a,0,0)\otimes(-x_a,0,0,0)]\\
=&[(-x_a,0,0,0)\otimes (1-b,0,0,0)+(1-b,0,0,0)\otimes (-x_a,0,0,0)]\\
=&[(x_a,0,0,0)\otimes (b-1,0,0,0)+(b-1,0,0,0)\otimes (x_a,0,0,0)].
\end{align*}
This corresponds to the element
\begin{equation}\label{eqn-ker-element-1}
  [(x_a,0)\otimes(b-1,0)+(b-1,0)\otimes(x_a,0)] \in \Z \otimes_{\Z\pi} \Gamma(M_1).
  \end{equation}
Here, similarly to above, we represent elements of $\Z \otimes_{\Z\pi} \Gamma(M_1)$ by elements of $(\Z\pi)^2 \otimes_{\Z} (\Z\pi)^2$.

If $m=2$, we will compute presently that $[2((1,0,0,0)\otimes(0,1,0,0)+(0,1,0,0)\otimes(1,0,0,0))]$ lies in the image of $\Z \otimes_{\Z\pi} \Gamma(M_1)$, so the order of $[(1,0,0,0)\otimes(0,1,0,0)+(0,1,0,0)\otimes(1,0,0,0)]$ in $\Z\otimes_{\Z\pi} E$ is at most 2, regardless of $n$.
To see this, we compute that
\begin{align*}
	&(2,0,0,0)\otimes(0,1,0,0)+(0,1,0,0)\otimes(2,0,0,0) \\
	= &(N_b+(1-b),0,0,0)\otimes(0,1,0,0)+(0,1,0,0)\otimes(N_b+(1-b),0,0,0)
\end{align*}
	and using the relations $(0,N_b,N_a,0)$ and $(N_a,0,0,0)$ in $\ker d_2$:
\begin{align*}
	 &[(N_b,0,0,0)\otimes(0,1,0,0)+(0,1,0,0)\otimes(N_b,0,0,0)] \\
   = &[(1,0,0,0)\otimes(0,N_b,0,0)+(0,N_b,0,0)\otimes(1,0,0,0)]\\
	=&[(1,0,0,0)\otimes(0,0,-N_a,0)+(0,0,-N_a,0)\otimes(1,0,0,0)]\\
	=&-[(N_a,0,0,0)\otimes(0,0,1,0)+(0,0,1,0)\otimes(N_a,0,0,0)] =0.
\end{align*}
Hence
\begin{align*}
	&[(2,0,0,0)\otimes(0,1,0,0)+(0,1,0,0)\otimes(2,0,0,0)] \\
	=&[(1-b,0,0,0)\otimes(0,1,0,0)+(0,1,0,0)\otimes(1-b,0,0,0)]\\
	=&[(1,0,0,0)\otimes(0,1-b,0,0)+(0,1-b,0,0)\otimes(1,0,0,0)].
\end{align*}
We computed above that this equals $[(1-b,0,0,0)\otimes(1-b,0,0,0)]$,
which corresponds to the element $[(1-b,0)\otimes(1-b,0)]\in \Z\otimes_{\Z\pi}\Gamma(M_1)$.
This completes the proof that $[2((1,0,0,0)\otimes(0,1,0,0)+(0,1,0,0)\otimes(1,0,0,0))]$ lies in the image of $\Z \otimes_{\Z\pi} \Gamma(M_1)$.

Similarly, with $x_b:=\sum_{i=1}^{m-1}(m-i)b^i$, we have
\begin{align*}
& [m((0,0,0,1)\otimes (0,0,1,0)+(0,0,1,0)\otimes(0,0,0,1))]\\
=& [(0,0,0, x_b) \otimes (0,0,0,a-1) + (0,0,0,a-1) \otimes (0,0,0,x_b)],
\end{align*}
which corresponds to the element
\begin{equation}\label{eqn-ker-element-2}
  [(0,-x_b)\otimes(0,(a-1)\lambda')+(0,(a-1)\lambda')\otimes(0,-x_b)]\in \Z \otimes_{\Z\pi} \Gamma(M_1).
\end{equation}
Moreover, if $n=2$, $2[(0,0,0,1)\otimes (0,0,1,0)+(0,0,1,0)\otimes(0,0,0,1)]$ is the image of $[(0,1-a)\otimes(0,1-a)]\in \Z\otimes_{\Z\pi}\Gamma(M_1)$.

Since the elements \eqref{eqn-ker-element-1} and \eqref{eqn-ker-element-2} (or their $m=2$, $n=2$ analogues) live in different direct summands in the decomposition
\[\Z\otimes_{\Z\pi}\Gamma(M_1)\cong \Z\otimes_{\Z\pi}\Gamma(\Z\pi/N_a)\oplus\Z\otimes_{\Z\pi}\Gamma(\Z\pi/N_b)\oplus (\Z\pi/N_a\otimes_{\Z\pi}\Z\pi/N_b),\]
we can check their divisibility separately and apply \cref{rem:torsion} to each summand.

 Again, the arguments are the same with the roles of $a$ and $b$ interchanged. Hence we only consider
\[[x_a\otimes(b-1)+(b-1)\otimes x_a]\in \Z\otimes_{\Z\pi}\Gamma(\Z\pi/N_a)\]
if $m\neq 2$ and $[(1-b)\otimes (1-b)]$ if $m=2$.

In the latter case we have $\Z\otimes_{\Z\pi}\Gamma(\Z \pi/N_a)\cong \Z\otimes_{\Z C_a}\Gamma(\Z C_a/N_a)\oplus \Z C_a/\langle N_a,g-\overline{g}\rangle$ by \cref{lem:torsion-free} with $G=C_b$, $H=C_a$, and $A= \Z C_a/N_a$.
This is a straightforward application of the lemma, once one has observed that
\[\Z \otimes_{\Z\pi} \big((\Z C_a/N_a \otimes_{\Z} \Z C_a /N_a)[C_b]/\textup{Flip}_b\big) \cong \Z C_a/\langle N_a,g-\overline{g}\rangle.\]
The image of $[(1-b)\otimes (1-b)]=2[1\otimes 1]+[1\otimes b+b\otimes 1]$ under the projection to $\Z C_a/\langle N_a,g-\overline{g}\rangle$ is $1$, which if $n$ is even is indivisible. To see that $1$ is indivisible, using the basis given by $[a^k]$ with $k\neq 0$, we have  $1=\sum_{i=1}^{n-1}a^i=\sum_{i=1}^{n/2-1}2a^i+a^{n/2}$. If $n$ is even, due to the coefficient of $a^{n/2}$ we see that $1$ is indivisible. If $n$ is odd, we obtain $1=\sum_{i=1}^{n-1}a^i=\sum_{i=1}^{(n-1)/2}2a^i$, so the element is divisible by~2.  However, for $n$ odd we have shown that both $n$ times and~2 times the image of the element  $[(1,0,0,0)\otimes(0,1,0,0)+(0,1,0,0)\otimes(1,0,0,0)]$ in $\Z \otimes_{\Z\pi} E$ is trivial. It follows that this image is in fact trivial in $\Z \otimes_{\Z\pi} E$, so we do not have to consider it for our application of \cref{lem:divisible}.

Now let us consider the case $m\neq 2$.	
Let $S$ be as in \cref{lem:torsion-free}, again with $G=C_b$, $H=C_a$, and $A= \Z[C_a]/N_a$. Assume that $b\in S$. By \cref{lem:torsion-free} we have
\[\Z\otimes_{\Z\pi}\Gamma(\Z\pi/N_a)\cong \Z\otimes_{\Z C_a}\Gamma(\Z C_a/N_a)\oplus\bigoplus_{S}\Z C_a/N_a\oplus L,\]
where $L\cong \Z C_a/\langle N_a,g-\overline{g}\rangle$ if the order $m$ of $b$ is even and $L=0$ otherwise.
Here we use the observation from above to find $L$ and that $\Z C_a/N_a \otimes_{\Z C_a} \Z C_a/N_a$.
In this direct sum decomposition, $[x_a\otimes(b-1)+(b-1)\otimes x_a]$ corresponds to
\[[x_a\otimes 1+1\otimes x_a]\in \Z\otimes_{\Z C_a}\Gamma(\Z C_a/N_a) \text{ plus }
x_a\in (\Z C_a/N_a)_{b\in S}.\]
In other words, the projection to $(\Z C_a/N_a)_{b\in S}$ is $x_a$.
Since $x_a\in \Z C_a/N_a$ is indivisible, so is $[x_a\otimes(b-1)+(b-1)\otimes x_a]$.
Now using \cref{lem:divisible}, we conclude that
there is no torsion in $\Z\otimes_{\Z\pi} \Gamma(\ker d_2)$.  By \cref{lemma:torsion-gamma-equals-tate} we therefore have $\wh{H}_0(\pi;\Gamma(\ker d_2))=0$.
\end{proof}

\subsection{The computation \texorpdfstring{of $\wh H_0(\pi;\Gamma(\coker d^2))$}{for the cokernel}}\label{sec:computation-cokernel}

In this section we show the following, which will complete the proof of \cref{thm:main}.

\begin{prop}
	\label{prop:coker}
	For every finite abelian group with two generators $\pi$, $\wh H_0(\pi;\Gamma(\coker d^2))=0$.
\end{prop}

Again we may assume that $|\pi|$ is a power of 2, but we also do not need this assumption for \cref{prop:coker}.

\begin{proof}
  Let $\pi=\langle a,b \mid a^{n}, b^{m}, [a,b] \rangle$.
As before, let $C_2\xrightarrow{d_2}C_1\xrightarrow{d_1}C_0$
be the chain complex corresponding to the presentation
$\langle a,b\mid a^n,b^m,[a,b]\rangle$.  Let $d^2 \colon C^1 = \Hom_{\Z\pi}(C_1,\Z\pi) \to C^2 = \Hom_{\Z\pi}(C_2,\Z\pi)$ be the dual of $d_2$.
The module
\[
	M \coloneqq \coker d^2
	\cong (\Z\pi)^3/\langle (N_a,1-b,0),(0,1-a,N_b) \rangle
\]
fits into an extension
\begin{equation}\label{eq:ses-M}
	0 \to \Z\pi/N\oplus\Z\pi/N\to M\to \Z \to 0,
\end{equation}
where $N$ is the norm element in $\Z\pi$.
The injection is given by the inclusion of the first and third summands in the presentation of $M$ above, given by $d^2$.
The surjection is given by the projection to the second summand.

By \cref{lem:bauer} there is a $\Z\pi$ module $D$ and short exact sequences
\begin{equation}\label{eq:M-ses-1}
  0 \to \Gamma(\Z\pi/N\oplus\Z\pi/N) \to \Gamma(M) \to D \to 0
\end{equation}
and
\begin{equation}\label{eq:M-ses-2}
0 \to (\Z\pi/N\oplus\Z\pi/N) \otimes_{\Z} \Z \to D \to \Gamma(\Z) \to 0.
\end{equation}
As in \cref{sec:computation-kernel}, we will show that the first satisfies
the conditions of \cref{lem:divisible},
and use the second to compute the torsion in $\Z \otimes_{\Z\pi} D$.
Then we will apply the first again, tensored with $\Z \otimes_{\Z \pi}$,
together with \cref{lem:divisible},
to deduce that $\Z \otimes_{\Z\pi} \Gamma(M)$ is torsion free.

\begin{claim}
  The short exact sequence~\eqref{eq:M-ses-1} satisfies the hypotheses of the first sentence of \cref{lem:divisible}. That is $\Gamma(M)$ and  $\Z \otimes_{\Z\pi} \Gamma(\Z\pi/N\oplus\Z\pi/N)$ are torsion free as $\Z$-modules.
\end{claim}

First, $\Z\pi/N\oplus\Z\pi/N$ and $\Z$ are torsion free abelian groups, and therefore so is $M$. Alternatively, apply \cref{lem:ker-and-coker-torsion-free}.
It follows that $\Gamma(M)$ is torsion free by \cref{lem:gammafree}.
Next, we have
\[
	\Gamma(\Z\pi/N\oplus\Z\pi/N)
	\cong \Gamma(\Z\pi/N)\oplus\Gamma(\Z\pi/N)\oplus (\Z\pi/N\otimes_\Z\Z\pi/N)
\]
By \cref{thm:torsion-free-gamma}~(iii),
$\Z \otimes_{\Z\pi} \Gamma(\Z\pi/N)$ is torsion free.
Additionally
\[
	\Z \otimes_{\Z\pi}(\Z\pi/N \otimes_{\Z} \Z\pi/N)
	\cong \Z\pi/N \otimes_{\Z\pi} \Z\pi/N
	\cong \Z\pi/N
\]
is torsion free (as an abelian group) and thus
$\Z\otimes_{\Z\pi}\Gamma(\Z\pi/N\oplus\Z\pi/N)$
is torsion free.
 This completes the proof of the claim that the short exact sequence~\eqref{eq:M-ses-1} satisfies the hypotheses of the first sentence of \cref{lem:divisible}.
 \newline

We therefore have a short exact sequence:
\begin{equation}\label{eq:ses-M-tensored-Z}
	0 \to \Z \otimes_{\Z\pi} \Gamma(\Z\pi/N\oplus\Z\pi/N)
	\to \Z \otimes_{\Z\pi}\Gamma(M)
	\to \Z \otimes_{\Z\pi} D
	\to 0.
\end{equation}
We will use this later in conjunction with \cref{lem:divisible},
to deduce that $\Z \otimes_{\Z\pi} \Gamma(M)$ is torsion free.

Next, we use the extension \eqref{eq:M-ses-2}
to compute the torsion in $\Z \otimes_{\Z\pi} D$.
Using $\Gamma(\Z) \cong \Z$, the extension gives a long exact sequence ending with:
\begin{equation}\label{eqn:Tate_boundary_map_application}
\cdots \to \widehat H_1(\pi;D) \to \widehat H_1(\pi;\Z) \xrightarrow{\partial}
  \widehat H_0(\pi;\Z\pi/N\oplus\Z\pi/N) \to
	 \widehat H_0(\pi;D) \to \widehat H_0(\pi;\Z) = 0.
\end{equation}
The above boundary map $\partial$ in the long exact sequence of Tate homology groups
is the boundary map from standard homology, noting that it has image in the
subgroup $\ker (\cdot N)$ of the orbits $A_{\pi} = H_{0}(\pi; A)$, where $A = \Z\pi/N\oplus\Z\pi/N$.
We compute $\partial$ as follows.  Let $C_*$ be a free $\Z\pi$-module resolution of $\Z$ with $C_1 \cong (\Z\pi)^2 \xrightarrow{(a-1,b-1)} C_0 \cong \Z\pi$. In the next diagram we show part of the short exact sequence of chain complexes from which one computes the map $\partial$.
\begin{center}
	\begin{tikzcd}
	C_{1} \otimes_{\Z \pi} (\Z\pi/N\oplus\Z\pi/N) \ar[r] \ar[d, "d_1 \otimes \id_{\Z\pi/N\oplus\Z\pi/N}"] &
	C_{1} \otimes_{\Z \pi} D \ar[r] \ar[d, "d_1 \otimes \id_{D}"] &
	C_{1} \otimes_{\Z \pi} \Z \ar[d, "d_1 \otimes \id_{\Z} =0"]\\
	C_{0} \otimes_{\Z \pi} (\Z\pi/N\oplus\Z\pi/N) \ar[r] &
	C_{0} \otimes_{\Z \pi} D \ar[r] &
	C_{0} \otimes_{\Z \pi} \Z
	\end{tikzcd}
\end{center}
We have
$	\wh H_1(\pi;\Z)
	\cong C_1 \otimes_{\Z\pi}\Z/\im (d_2\otimes\id_{\Z})
	\cong \Z/n\oplus \Z/m $
and
$	\widehat H_0(\pi;\Z\pi/N\oplus\Z\pi/N)
	\cong \Z/|\pi|\oplus \Z/|\pi|. $
In the extension \eqref{eq:M-ses-2}, a preimage of $1\in \Z \cong \Gamma(\Z)$ in
$D \cong \frac{\Gamma(M)}{\Gamma(\Z\pi/N \oplus 0 \oplus \Z\pi/N)}$
is given by
$[(0,1,0)\otimes(0,1,0)]$.
Here we represent elements of $\Gamma(M)$ by elements of $(\Z\pi)^3 \otimes_{\Z} (\Z\pi)^3$; these in turn determine elements of $D$ since $D$ is a quotient of $\Gamma(M)$.
Under the map
\begin{center}
	\begin{tikzcd}
	C_{1} \otimes_{\Z \pi} D \ar[d, "d_1 \otimes \id_{D}"] \ar[r, equal] &
	(\Z\pi)^2 \otimes_{\Z\pi} D \cong D \oplus D &
	{([-(0,1,0)\otimes(0,1,0)], 0)} \ar[d, mapsto] \\
	C_{0} \otimes_{\Z \pi} D \ar[r, equal] &
	\Z\pi \otimes_{\Z\pi} D \cong D &
	{[(a-1)((0,1,0) \otimes (0,1,0))]}
	\end{tikzcd}
\end{center}
the element $[-(0,1,0)\otimes(0,1,0)]$ in the first factor of $D \oplus D$
is sent by $d_1 \otimes \id_D$ to $[(a-1)((0,1,0)\otimes(0,1,0))]$.
This tells us the image of $(-1,0) \in \Z/n \oplus 0 \subseteq \Z/n \oplus \Z/m \cong \wh H_1(\pi;\Z)$ under $\partial$. To describe this image precisely, we consider the following diagram.
\begin{center}
  \begin{tikzcd}
  & & \Z\pi/N\oplus\Z\pi/N \ar[d,"p"] \ar[r,hookrightarrow, "\iota"] & D \\
    \wh H_1(\pi;\Z) \ar[r,"\partial"] & \wh{H}_0 (\pi;\Z\pi/N \oplus \Z\pi/N) \ar[r,hookrightarrow,"T"] & \Z \otimes_{\Z\pi} (\Z\pi/N \oplus \Z\pi/N) &
    \end{tikzcd}
\end{center}
Here, the map $\iota$ is the identification  $\Z\pi/N\oplus\Z\pi/N \cong (\Z\pi/N\oplus\Z\pi/N)\otimes_{\Z} \Z$ followed by the map from \eqref{eq:M2-ses-2}, and $T$ is the inclusion of the torsion subgroup by \cref{lemma:torsion-gamma-equals-tate}.
Then the computations above show that $(-1,0) \in \Z/n \oplus 0 \subseteq \Z/n \oplus \Z/m \cong \wh H_1(\pi;\Z)$
is sent under  the composition $T \circ \partial$ to \[p(\iota^{-1}((a-1)((0,1,0)\otimes(0,1,0)))).\]
Here as above the element $(a-1)((0,1,0)\otimes(0,1,0)) \in D$ is represented as an element in $\Gamma(M)$, which in turn is represented by an element of $(\Z\pi)^3 \otimes_{\Z} (\Z\pi)^3$.
%
%
%
To express this as an element in
$\Z \otimes_{\Z\pi} (\Z\pi/N\oplus\Z\pi/N)$, we will now compute.
Using the relation $[(0,1-a,N_b)]=0$ in $M$ we have
\begin{align*}
0=& [(0,1-a,N_b) \otimes (0,-a,0) - (0,1,N_b)\otimes (0,1-a,N_b) ]\\
 =& [(0,1-a,N_b)\otimes(0,-a,0)-(0,1,N_b)\otimes(0,-a,0)-(0,1,N_b)\otimes(0,1,N_b)]\\
 =& [(0,1,N_b)\otimes(0,-a,0) + (0,-a,0)\otimes(0,-a,0) -(0,1,N_b)\otimes(0,-a,0)-(0,1,N_b)\otimes(0,1,N_b)]\\
 =& [(0,-a,0)\otimes(0,-a,0)-(0,1,N_b)\otimes(0,1,N_b)]\\
 =& [(0,a,0)\otimes(0,a,0)-(0,1,0)\otimes(0,1,0)-(0,0,N_b)\otimes(0,0,N_b)\\
  & \; -(0,1,0)\otimes(0,0,N_b)-(0,0,N_b)\otimes(0,1,0)].
\end{align*}
Hence we have
\begin{align*}
	 [(a-1)((0,1,0)\otimes(0,1,0))]
   = [(0,0,N_b)\otimes(0,0,N_b)+(0,1,0)\otimes(0,0,N_b)+(0,0,N_b)\otimes(0,1,0)].
\end{align*}
As $D$ is the quotient of
$\Gamma(M) = \Gamma(\frac{\Z\pi \oplus \Z\pi \oplus \Z\pi}{\langle \cdots \rangle})$
by $\Gamma(\Z\pi/N \oplus 0 \oplus \Z\pi/N)$, we have in $D$:
\begin{align*}
	  [(a-1)((0,1,0)\otimes(0,1,0))]
	= [(0,1,0)\otimes(0,0,N_b) + (0,0,N_b)\otimes(0,1,0)].
\end{align*}
This element has preimage $(0,N_b)\otimes 1\in (\Z\pi/N\oplus\Z\pi/N)\otimes_\Z \Z$ and thus when tensoring with $\Z$ over $\Z\pi$, we get
$(0,m)\in \Z/|\pi|\oplus \Z/|\pi| \cong \widehat H_0(\pi;\Z\pi/N\oplus\Z\pi/N) \subseteq \Z \otimes_{\Z\pi} (\Z\pi/N\oplus\Z\pi/N)$.

Similarly, the generator of $\Z/m$ in $\Z/n \oplus \Z/m \cong \wh H_1(\pi;\Z)$
maps under the boundary map $\partial$ to
$(n,0)\in \Z/|\pi|\oplus \Z/|\pi|$.
Hence the torsion in $\Z \otimes_{\Z\pi} D$, which equals $\widehat H_0(\pi;D) \cong \coker \partial$ by \cref{lemma:torsion-gamma-equals-tate} and \eqref{eqn:Tate_boundary_map_application}, is isomorphic to $\Z/n \oplus \Z/m$.

Now recall the short exact sequence~\eqref{eq:ses-M-tensored-Z},
	$0 \to \Z \otimes_{\Z\pi} \Gamma(\Z\pi/N \oplus \Z\pi/N)
	\xrightarrow{\psi} \Z \otimes_{\Z\pi} \Gamma(M)
	\xrightarrow{\theta} \Z \otimes_{\Z\pi} D
	\to 0$.
We want to apply the second half of \cref{lem:divisible}
to deduce that $\Z\otimes_{\Z\pi}\Gamma(M)$ is torsion free.
The elements $u:= [(0,1,0) \otimes (0,0,1) + (0,0,1) \otimes (0,1,0)]$
and $v:=[(0,1,0) \otimes (1,0,0) + (1,0,0) \otimes (0,1,0)]$
are preimages under $\theta$
in $\Z\otimes_{\Z\pi}\Gamma(M)$ of the generators
of the torsion in $\Z\otimes_{\Z\pi} D$.

Next we find the preimage under $\psi$ of $m$ times $u=[(0,1,0) \otimes (0,0,1) + (0,0,1) \otimes (0,1,0)]$.
The first equation of the following
computation in $\Z \otimes_{\Z \pi} \Gamma(M)$ uses
the identity $m-N_b=-y_b(1-b)$
in $\Z \pi$ with
$y_b \coloneqq \sum_{i=0}^{m-1} i b^{i} = b + 2b^{2} + \cdots + (m-1)b^{m-1}$.
\begin{align*}
	& [m ((0,1,0)\otimes(0,0,1)+(0,0,1)\otimes (0,1,0))] \\
	=& [(0,N_b-y_b(1-b),0)\otimes(0,0,1)+(0,0,1)\otimes (0,N_b-y_b(1-b),0)] \\
	=& [(N_ax_b,N_b,0)\otimes(0,0,1)+(0,0,1)\otimes (N_ax_b,N_b,0)],
\end{align*}
where the last equation uses that $0=[(N_a,1-b,0)]$ in $M$.
Since $0=[(a-1)((0,1,0)\otimes(0,1,0))] \in \Z\otimes_{\Z\pi}\Gamma(M)$,
the above computation for the boundary map $\partial$ yields
\[
	[-(0,0,N_b)\otimes(0,0,N_b)]
	=[ (0,1,0)\otimes(0,0,N_b) + (0,0,N_b)\otimes(0,1,0)].
\]
Since we tensored with $\Z$ over $\Z\pi$, where $\Z\pi$ acts diagonally on $(\Z\pi)^3 \otimes_{\Z} (\Z\pi)^3$,
we have
\begin{align*}
	 & [(0,1,0)\otimes(0,0,N_b)+(0,0,N_b)\otimes(0,1,0)]\\
= &[(0,\overline{N}_b,0)\otimes(0,0,1)+(0,0,1)\otimes(0,\overline{N}_b,0)]\\
= &[(0,N_b,0)\otimes(0,0,1)+(0,0,1)\otimes(0,N_b,0)].
\end{align*}
Combining these we get
\begin{align*}
	& [m ((0,1,0)\otimes(0,0,1)+(0,0,1)\otimes (0,1,0))] \\
	=& [(N_ax_b,0,0)\otimes(0,0,1)+(0,0,1)\otimes(N_ax_b,0,0)-(0,0,N_b)\otimes(0,0,N_b)].
\end{align*}

To show that an element of a direct sum is indivisible,
it is enough to check this in one summand.
Hence consider the image
under the projection
\begin{center}
	\begin{tikzcd}
	\Z\otimes_{\Z\pi}\Gamma(\Z\pi/N\oplus\Z\pi/N) \ar[r, "\cong"] &
	\Z \otimes_{\Z \pi}
	(\Gamma(\Z\pi/N) \oplus \Gamma(\Z\pi/N) \oplus (\Z\pi/N \otimes_{\Z} \Z\pi/N))
	\ar[d, twoheadrightarrow]\\
	& \Z \otimes_{\Z \pi} (\Z\pi/N \otimes_{\Z} \Z\pi/N) \cong \Z\pi/N
	\end{tikzcd}
\end{center}
which for the preimage of $mu$ is $N_a y_b$.  An analogous computation shows that the preimage under $\psi$ of $n$ times $v=[(0,1,0) \otimes (1,0,0) + (1,0,0) \otimes (0,1,0)]$ projects to $N_b y_a \in \Z\pi/N \cong \Z \otimes_{\Z \pi} (\Z\pi/N \otimes_{\Z} \Z\pi/N)$, where $y_a \coloneqq \sum_{i=0}^{n-1} i a^{i}$.

\begin{claim}
To complete the proof it suffices to show, for every $k,\ell$ with $\gcd(k,\ell)=1$, that $k N_a y_b+\ell N_b y_a$ is indivisible.
\end{claim}
To see this, let $a,b \in \Z$ and write $K:= \gcd(na,mb)$.
Then to check that $au+bv$ is not torsion in $\Z \otimes_{\Z\pi} \Gamma(M)$, as per \cref{lem:divisible}, note that $mn/K$ is the order of the image of $au+bv$ in $\Z \otimes_{\Z\pi} D$. Indeed, we have
\[\left(\frac{mn}{K}\right) (au+bv) = \left(\frac{na}{K}\right)mu +  \left(\frac{mb}{K}\right)nv,\]
 and this maps to $0$ in $\Z \otimes_{\Z\pi} D$ since both $mu$ and $nv$ do.  Taking the preimage in $\Z\otimes_{\Z\pi}\Gamma(\Z\pi/N\oplus\Z\pi/N)$, and then projecting to $\Z \otimes_{\Z \pi} (\Z\pi/N \otimes_{\Z} \Z\pi/N) \cong \Z\pi/N$, we obtain
\[\left(\frac{na}{K}\right) N_a y_b +  \left(\frac{mb}{K}\right) N_b y_a.\]
Write $k:= na/K$ and $\ell := mb/K$. Note that $\gcd(k,\ell)=1$, since $K= \gcd(na,mb)$.
Then it suffices to show that $k N_a y_b + \ell N_b y_a$ is indivisible.  In particular, we see that it is enough to show that $k N_a y_b + \ell N_b y_a$ is indivisible for every pair of coprime integers $k,\ell$. This completes the proof of the claim.

Now we show the desired statement.
As an abelian group, $\Z\pi/N$ is isomorphic to the submodule
of elements in $\Z\pi$ whose coefficient at the identity group element is zero.
Observe that $N_a y_b$ and $N_b y_a$ lie in this submodule.
Let $k,\ell \in \Z$ be coprime integers.
  In $kN_a y_b+\ell N_b y_a$, the coefficient of $a^1 b^0 = a$ is $\ell$,
	while the coefficient of $a^0 b^1 = b$ is $k$
	(and the coefficient in front of the trivial group element is zero).  It follows that
	in $\Z\pi/N$, $kN_a y_b+\ell N_b y_a$ can only be divisible by
	$\gcd(k,\ell)$ and its divisors. But $\gcd(k,\ell)=1$ so $kN_a y_b+\ell N_b y_a$ is indivisible, as desired.
Hence there is no torsion in $\Z\otimes_{\Z\pi}\Gamma(M)$
by \cref{lem:divisible}. Therefore $\wh H_0(\pi;\Gamma(M)) = \wh H_0(\pi;\Gamma(\coker d^2))=0$ by \cref{lemma:torsion-gamma-equals-tate}.
\end{proof}

\section{\texorpdfstring{Groups of order at most $16$}{Groups of order at most 16}}\label{sebsection:groups-less-16}
For small groups the following computations were made
by Hennes \cite{hennes} using a computer program.
For $\pi=(\Z/2)^2,\Z/2\times\Z/4,(\Z/4)^2$ and $D_8$, $\Z\otimes_{\Z\pi}\Gamma(\ker d_2\oplus \coker d^2)$ is torsion free.
For $\pi=(\Z/2)^3, \Z/2\times D_8$ and
$\Z/4\times \Z/2\times \Z/2$, $\Z\otimes_{\Z\pi}\Gamma(\ker d_2\oplus \coker d^2)$
contains torsion.

For $\pi=(\Z/2)^2$ it was claimed in \cite[Final~remark]{bauer} (without proof),
that $\wh H_0(\pi;\Gamma(\coker d^2))\cong (\Z/2)^2$. By Hennes' computation and also by our proof above, this is not true.

We used an algorithm written by the third author
in his bachelor thesis \cite{ruppik} to calculate
$\Z\otimes_{\Z\pi}\Gamma(\ker d_2)$ for all groups $\pi$ up to order 16.
The algorithm was written in SageMath
and the source code together with output from
the calculations
can be found in the GitHub repository \cite{algorithm}.
There were no computations made of $\Z\otimes_{\Z\pi}\Gamma(\coker d^2)$.

\begin{example}
	The following is a complete list of all groups of order at most $16$
	such that $\wh H_0(\pi;\Gamma(\ker d_2))$
	is nontrivial.
	The group $Q_8=\langle i,j,k\mid i^2=j^2=k^2=ijk\rangle$
	is the quaternion group.
	\begin{center}
		\begin{tabular}{l|l}
			$\pi$&$\wh H_0(\pi;\Gamma(\ker d_2))$  \\
			\hline
			\rule{0pt}{\normalbaselineskip}$\Z/4\times \Z/2\times \Z/2$&$(\Z/2)^2$ \\
			$\Z/2\times \Z/2\times \Z/2\times \Z/2$&$(\Z/2)^4$\\
			$Q_8\times \Z/2$&$(\Z/2)^4$
		\end{tabular}
	\end{center}
\end{example}

\bibliographystyle{alpha}
\bibliography{gamma}
\end{document}